\documentclass[11pt, oneside]{article}   	
\usepackage{amsmath,amssymb,amsthm,geometry,hyperref}
\usepackage{hyperref}
\usepackage{graphicx}						
\usepackage{amssymb}
\usepackage[utf8]{inputenc}
\usepackage[T1]{fontenc}
\usepackage{listings}
\usepackage{xcolor}
\usepackage{diagbox}
\usepackage{placeins}
\usepackage{subfig}
\usepackage{multirow}
\usepackage[justification=centering]{caption}

\newcommand{\be}{\begin{equation}}
\newcommand{\ee}{\end{equation}}
\newcommand{\bea}{\begin{eqnarray}}
\newcommand{\eea}{\end{eqnarray}}
\newcommand{\bean}{\begin{eqnarray*}}
\newcommand{\eean}{\end{eqnarray*}}
\newcommand{\nn}{\nonumber}

\newcommand{\si}{\sigma}

\newcommand{\pa}{\partial}

\newcommand{\footremember}[2]{%
    \footnote{#2}
    \newcounter{#1}
    \setcounter{#1}{\value{footnote}}%
}
\newcommand{\footrecall}[1]{%
    \footnotemark[\value{#1}]%
} 

\let\Oldsection\section
\renewcommand{\section}{\FloatBarrier\Oldsection}

\let\Oldsubsection\subsection
\renewcommand{\subsection}{\FloatBarrier\Oldsubsection}

\let\Oldsubsubsection\subsubsection
\renewcommand{\subsubsection}{\FloatBarrier\Oldsubsubsection}

\providecommand{\keywords}[1]{\textbf{\textit{Keywords---}} #1}

\title{Revisiting Feller Diffusion: Derivation and Simulation}
\author{%
  Ranjiva Munasinghe\footremember{MINDAM}{MIND Analytics \& Management, 10/1 De Fonseka Place, Colombo 5, Sri Lanka.}%
  \and Leslie Kanthan\footremember{UCL}{University College London, Gower Street, London, WC1E 6BT, UK.} \footnote{Data Spartan, 88 Wood Street, St. Pauls, London, EC2V 7RS, UK.}%
  \and Pathum Kossinna\footrecall{MINDAM} %
  }					

\begin{document}

\maketitle	

\begin{abstract}
We propose a simpler derivation of the probability density function of Feller Diffusion using the Fourier Transform and solving the resulting equation via the Method of Characteristics. We also discuss simulation algorithms and confirm key properties related to hitting time probabilities via the simulation.
\end{abstract}

\keywords{Feller Diffusion, Bessel Squared Process of Dimension Zero, Zero Degrees of Freedom, Square Root Diffusion, Simulation, Method of Characteristics, Absorbing Random Walks, Survival Probability, Hitting Time, First Passage Time.}

\section{Introduction}
Feller Diffusion, also known as Bessel Squared Process of dimension zero \cite{Feller, Feller2, Feller3, Wax, GL}, is given by the \^{I}to stochastic differential equation (SDE):
\bea
\dot{x}(t) &=& \si \sqrt{x} \cdot \xi(t) \label{eqn:SDE}\\
x(t_o) &=& x_o > 0 \label{eqn:IC} 
\eea
Here $x(t)$ denotes the position of a random walker following the trajectory given by equation (\ref{eqn:SDE}) with initial condition $x_o$ at (initial) time $t_o$, and $\xi(t)$ is a Gaussian white noise independent of $x(t)$, with amplitude parameter $\si $, and the following correlation functions describing its statistics:
\bea
\langle \xi(t) \rangle_{\xi} &=& 0 \\
\langle \xi(t) \xi(t')\rangle_{\xi} &=& \delta(t-t') 
\eea
The averaging $\langle \cdots \rangle_{\xi}$ denotes an average with respect to the realizations of the stochastic process $\xi(t)$. Due to presence of the square root term the dynamics are such that the random walker's postion is always positive and if it reaches the origin it will stay there (absorbing barrier). Feller Diffusion is in the class of \^{I}to SDEs known as square root diffusions \cite{Glass} and is a special case in that it is driftless.  In general, square root diffusions have a drift term in addition to the multiplicative square root noise term. Feller's motivation in this area was genetics \cite{Feller2} with natural extensions to biological systems and population dynamics. Square root diffusions also have application in physics - in particular in the area of branching processes and percolation theory \cite{Cardy} and phase transitions with absorbing states \cite{HH}. In addition to the science applications listed above, there are applications of square root diffusions in quantitative finance - e.g. the CEV Model for asset pricing \cite{Glass} and CIR Model for interest rate modeling \cite{Hull, Wilmott}.  
\\
The Fokker-Planck equation corresponding to the SDE (\ref{eqn:SDE}) and initial condition (\ref{eqn:IC}) is given by
\bea
\frac{\pa}{\pa t} P(x,t) &=& \frac{\si^2}{2} \frac{\pa^2}{\pa x^2} [ x P(x,t) ] \label{eqn:FPE} \\
P(x,t_o) &=& \delta(x-x_o) \\
\int_{0}^{\infty} P(x,t) \; dx &=& 1.
\eea
Here $P(x,t) = P(x,t | x_o, t_o) $ describes the probability density of the position of the random walker at a particular time via the following relation:
\[
P(x,t)  \Delta x = \text{Probability walker} \in  [x,x+\Delta x]  \text{ at time }  t
\]

\section{Derivation of the Probability Density}
The first step in deriving the probability density is to use the Fourier transform on $P(x,t)$:
\[
\hat{P}(k,t) = \int_{-\infty}^{\infty} P(x,t) e^{-ikx} dx = \int_{0}^{\infty} P(x,t) e^{-ikx} dx
\]
The Fokker-Planck equation (\ref{eqn:FPE}) now becomes
\bea
\frac{\pa}{\pa t} \hat{P}(k,t) &=& \frac{-i \si^2 k^2}{2} \frac{\pa}{\pa k} \hat{P}(k,t) \label{eqn:kFPE} \\
\hat{P}(k,t_o) &=& \exp \left[ -ikx_o \right] \label{eqn:kIC} \\
\hat{P}(k=0,t)  &=& 1. \label{eqn:kNorm}
\eea
This partial differential equation (PDE) can be solved by the method of characteristics (please refer Appendix \ref{appx:MoC}) to yield
\be
\label{eqn:soln}
\hat{P}(k,t) = \exp \left[ \frac{-ikx_o}{1 + \frac{ik \si^2}{2} (t-t_o)} \right]
\ee
Note that the initial condition given by equation (\ref{eqn:kIC}) and normalization given by equation (\ref{eqn:kNorm}) are easily verified.
The characteristic function $\varphi(k,t)$ can be easily derived by the transformation $k \rightarrow -k$ in the Fourier transform, i.e. using the relation:
\[
\varphi(k,t) = \hat{P}(-k,t)
\]
This gives the characteristic function for the probability density governing the SDE (\ref{eqn:SDE}) with initial condition (\ref{eqn:IC}) as
\be
\label{eqn:char}
\varphi(k,t) = \exp \left[ \frac{ikx_o}{1 - \frac{ik \si^2}{2} (t-t_o)} \right]
\ee
This is the characteristic function of a non-central chi-squared distribution with zero degrees of freedom (please refer to Appendix \ref{appx:0DNCXSq}) where 
\begin{enumerate}

\item The random variable has been scaled such that $x \rightarrow \frac{4}{\si^2 (t-t_o)} x$ \label{item:xscale}

\item The non-centrality parameter is $\lambda = \frac{4x_o}{\si^2 (t-t_o)}$ \label{item:noncen}

\item Then $P(x,t | x_o, t_o ) =  g \left( \frac{4x}{\si^2 (t-t_o)}; \frac{4x_o}{\si^2 (t-t_o)} \right)$, where $g(x;\lambda)$ is the density function of a non-central chi-squared distribution with zero degrees of freedom with non-centrality parameter $\lambda$. 

\end{enumerate}
The non-centrality parameter \cite{Siegel} enables us to estimate the probability that the random variable takes the value 0 as $e^{-\lambda /2}$. In terms of the stochastic process (\ref{eqn:SDE}) this is the probability the random walker gets absorbed \cite{GL}  at the origin at time $t$ given the starting point $x_o$. We express these in terms of the survival function\footnote{$S(t)$ is the probability the random walker has survived up to time $t$.}, $S(t)$ and the corresponding absorption probability\footnote{$F(t) = 1 - S(t)$ and is the probability the random walker has been absorbed by time $t$.} $F(t)$:
\be
\label{eqn:absprob}
F(t) = \exp \left[ \frac{ -2x_o}{\si^2 (t-t_o)} \right]
\ee
Note that in the limit $t \rightarrow \infty$ this probability will tend to 1 which means that all random walkers following process (\ref{eqn:SDE}) will eventually be absorbed\footnote{Alternatively $\lim_{t \rightarrow \infty} S(t) \rightarrow 0$.} at the origin.

\section{Sample Path Simulation}
\begin{figure}[h]
	\centering
	\includegraphics[width=0.75\textwidth]{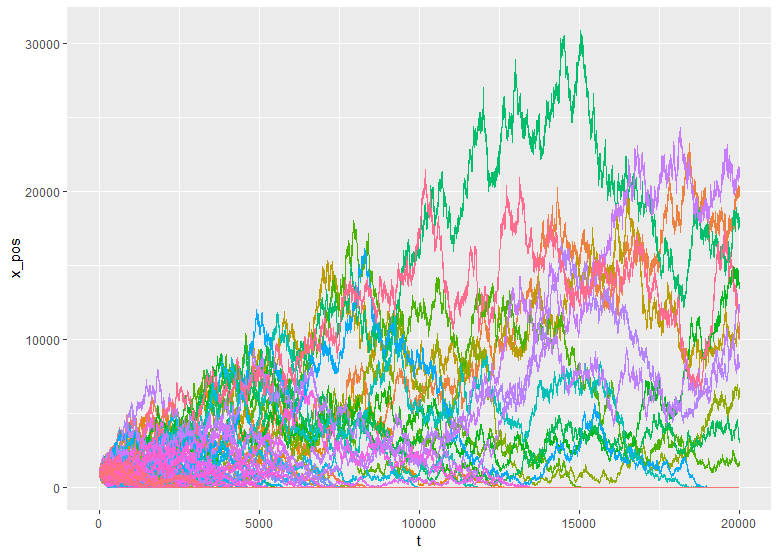}
	\caption{Feller Sample Path Simulation with $t_0=0$, $t_n=20000$, $\Delta t_i=1$, $\sigma^2=1$ and $x_0=1000$ (100 paths)}
	\label{fig:samplesim1}
\end{figure}

Using the fact that a non-central chi-squared distribution can be  expressed as a Poisson compound mixture of central chi-square distribution (refer Appendix \ref{appx:0DNCXSq}) with even degrees of freedom \cite{Siegel} leads to a very simple formulation for the simulation of a random walk \cite{Glass} corresponding to (\ref{eqn:SDE}). We make use of items \ref{item:xscale} and \ref{item:noncen} along with the Markov property \cite{GL, Oks} of the random walkers to describe the algorithm for generating a sample path:
\begin{enumerate}

\item For $i=1,2,\ldots,n$: 
\bean
c &\leftarrow& \frac{\si^2}{4 \Delta t_i} \;\;\; \text{where } \Delta t_i = t_i - t_{i-1} \\
\lambda &\leftarrow& x_{i-1} / c
\eean

\item Generate $N \sim$ Poisson ($\lambda / 2$) 

\item If $N=0$ set $X=0$ otherwise generate $X \sim \chi^2_{2N}$ 

\item Set $x_i = cX$ 

\end{enumerate}
The algorithm illustrated above was implemented in Rcpp (an implementation of C++ in R) and R was used in generating plots and data manipulation. The Rcpp implementation is given in Appendix \ref{appx:RCPP}. A simple simulation of 100 such generated paths is shown in Figure \ref{fig:samplesim1}.

\subsection{Examining Survival/Absorption Probabilities}
Using our simulations we can explore several properties of Feller Diffusion paths\footnote{Henceforth we set $t_o=0$ in both our simulation and computation.}. We begin by looking at the absorption probability (\ref{eqn:absprob}) and estimate this via the evolution of the (number of) simulated surviving paths. Our simulations match the intuition - i.e. more volatile walkers (i.e. larger values of $\si$) have lower survival rate and for large times the survival drops. Further the values are in close agreement with the theoretical values (Table \ref{tab:samplesim1} - using 10,000 simulations). 

\begin{table}[h]
	\centering
	\small
	\begin{tabular}{l|r|r|r|r|r|r|r|r|r|r|}
		\cline{2-11}
		& \multicolumn{10}{c|}{\textbf{Time (t)}}                                                                                                                                                                                                                                                                          \\ \hline
		\multicolumn{1}{|c|}{\multirow{2}{*}{\textbf{$\si^2$}}} & \multicolumn{2}{c|}{100}                                  & \multicolumn{2}{c|}{1,000}                                & \multicolumn{2}{c|}{5,000}                                & \multicolumn{2}{c|}{10,000}                               & \multicolumn{2}{c|}{20,000}                               \\ \cline{2-11} 
		\multicolumn{1}{|c|}{}                              & \multicolumn{1}{c|}{Theo.} & \multicolumn{1}{c|}{Sim.} & \multicolumn{1}{c|}{Theo.} & \multicolumn{1}{c|}{Sim.} & \multicolumn{1}{c|}{Theo.} & \multicolumn{1}{c|}{Sim.} & \multicolumn{1}{c|}{Theo.} & \multicolumn{1}{c|}{Sim.} & \multicolumn{1}{c|}{Theo.} & \multicolumn{1}{c|}{Sim.} \\ \hline
		\multicolumn{1}{|l|}{0.1}                           & 0.00                      & 0.00                          & 0.00                      & 0.00                          & 1.83                      & 1.70                          & 13.53                     & 13.47                         & 36.79                     & 36.11                         \\
		\multicolumn{1}{|l|}{1}                             & 0.00                      & 0.00                          & 13.53                     & 13.31                         & 67.03                     & 67.08                         & 81.87                     & 82.10                         & 90.48                     & 90.68                         \\
		\multicolumn{1}{|l|}{10}                            & 13.53                     & 13.02                         & 81.87                     & 81.84                         & 96.08                     & 96.12                         & 98.02                     & 98.03                         & 99.00                     & 99.03                         \\
		\multicolumn{1}{|l|}{100}                           & 81.87                     & 82.10                         & 98.02                     & 97.97                         & 99.60                     & 99.61                         & 99.80                     & 99.76                         & 99.90                     & 99.88                         \\ \hline
	\end{tabular}
	\caption{Percentage of paths absorbed at origin over time for different values of $\si^2$ with $t_0=0$, $\Delta t_i=1$, $x_0=1000$. (10,000 paths)}
	\label{tab:samplesim1}
\end{table}
We also observe that the further away from the origin the starting point (larger values of $x_o$) the better the chance of survival - as per the intuition from equation (\ref{eqn:absprob}). These results are displayed in  (Table \ref{tab:samplesim4} - using 10,000 simulations) and show close agreement with the theoretical values:
\begin{table}[h]
	\centering
	\small
	\begin{tabular}{l|r|r|r|r|r|r|r|r|r|r|}
		\cline{2-11}
		& \multicolumn{10}{c|}{Time (t)}                                                                                                                                                                                                                                                             \\ \hline
		\multicolumn{1}{|c|}{\multirow{2}{*}{$x_0$}} & \multicolumn{2}{c|}{100}                               & \multicolumn{2}{c|}{1,000}                             & \multicolumn{2}{c|}{5,000}                             & \multicolumn{2}{c|}{10,000}                            & \multicolumn{2}{c|}{20,000}                            \\ \cline{2-11} 
		\multicolumn{1}{|c|}{}                       & \multicolumn{1}{c|}{Theo.} & \multicolumn{1}{c|}{Sim.} & \multicolumn{1}{c|}{Theo.} & \multicolumn{1}{c|}{Sim.} & \multicolumn{1}{c|}{Theo.} & \multicolumn{1}{c|}{Sim.} & \multicolumn{1}{c|}{Theo.} & \multicolumn{1}{c|}{Sim.} & \multicolumn{1}{c|}{Theo.} & \multicolumn{1}{c|}{Sim.} \\ \hline
		\multicolumn{1}{|l|}{10}                     & 81.87                    & 82.10                   & 98.02                    & 97.97                   & 99.60                    & 99.61                   & 99.80                    & 99.76                   & 99.90                    & 99.87                   \\ \hline
		\multicolumn{1}{|l|}{100}                    & 13.53                    & 13.02                   & 81.87                    & 81.85                   & 96.08                    & 96.11                   & 98.02                    & 98.04                   & 99.00                    & 99.03                   \\ \hline
		\multicolumn{1}{|l|}{1000}                   & 0.00                     & 0.00                    & 13.53                    & 13.31                   & 67.03                    & 67.08                   & 81.87                    & 82.10                   & 90.48                    & 90.68                   \\ \hline
		\multicolumn{1}{|l|}{10000}                  & 0.00                     & 0.00                    & 0.00                     & 0.00                    & 1.83                     & 1.79                    & 13.53                    & 13.51                   & 36.79                    & 36.07                   \\ \hline
	\end{tabular}
	\caption{Percentage of paths absorbed at origin over time for different values of $x_0$ with $t_0=0$, $\si^2=1$, $\Delta t_i=1$. (10,000 paths)}
	\label{tab:samplesim4}
\end{table}

\begin{figure}[h]
	\centering
	\subfloat[Feller Sample Path Simulation with $\sigma^2=0.1$]{\label{fig:a}\includegraphics[width=0.45\linewidth]{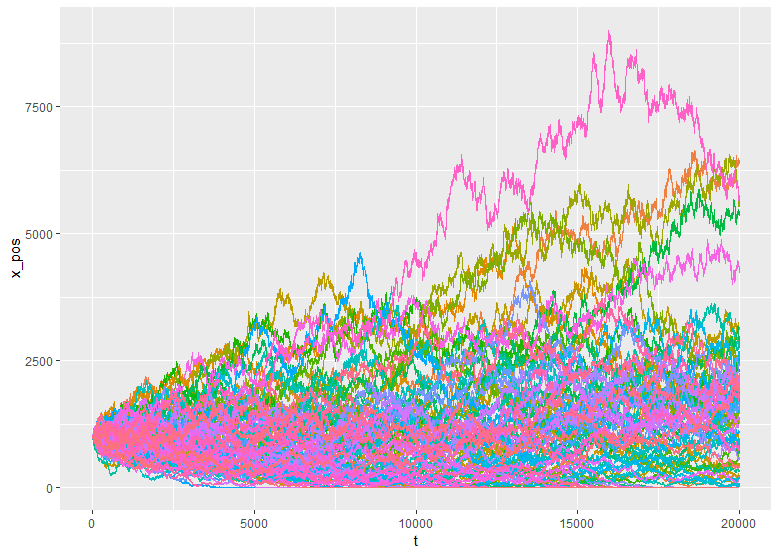}}\qquad
	\subfloat[Feller Sample Path Simulation with $\sigma^2=1$]{\label{fig:b}\includegraphics[width=0.45\linewidth]{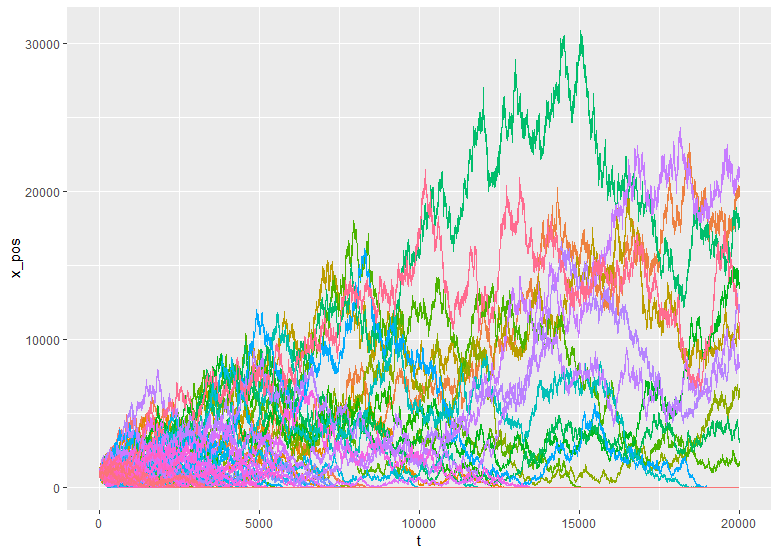}}\\
	\subfloat[Feller Sample Path Simulation with $\sigma^2=10$]{\label{fig:c}\includegraphics[width=0.45\textwidth]{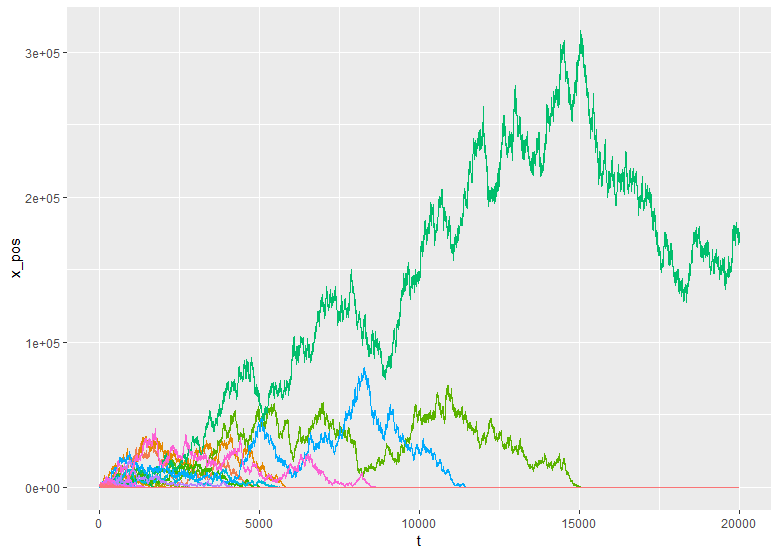}}\qquad%
	\subfloat[Feller Sample Path Simulation with $\sigma^2=100$]{\label{fig:d}\includegraphics[width=0.45\textwidth]{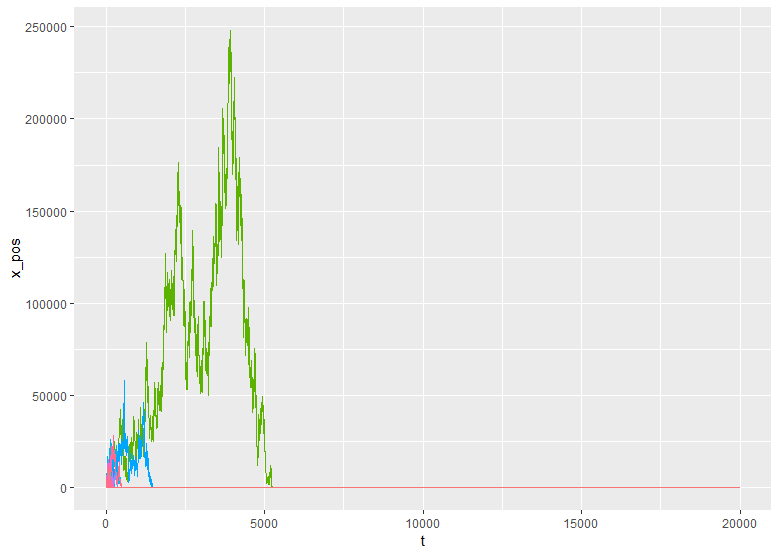}}%
	\caption{Feller Sample Path Simulations with $t_0=0$, $t_n=20000$, $\Delta t_i=1$, $x_0=1000$ and different values of $\si^2$ (100 paths).}
	\label{fig:varyingSigma}
\end{figure}

\subsection{Martingale Property and Average Position} \label{sec:MartingaleProperty}
As equation (\ref{eqn:SDE}) is a martingale\footnote{An SDE with no drift term.}, the expected position\footnote{The averaging is with respect to $P(x,t)$.} of a walker at time $t$ to be $x_o$ - calculation can be found in the Appendix \ref{appx:MeanVarFDX}. We can examine this via our simulation by evaluating the average position of several walkers at a particular time via:
\be
\label{eqn:simmeanpos}
\overline{x} (t) = \frac{1}{N} \sum_{m=1}^{N} x_m (t)
\ee
The fact that $\langle x(t) \rangle = x_o$ (given $x(0) = x_o$) might seem paradoxical given that the Feller process is one where eventually all walkers are absorbed at the origin. This can be understood through the simulation where the walkers that survive (even though this number will be small) will most likely be very far away from the origin and thus driving up the average given by equation (\ref{eqn:simmeanpos}). Results are displayed in Tables \ref{tab:avgxvar10} and \ref{tab:avgxvar1}. Note the error bounds (confidence intervals) for the Monte Carlo estimates are in agreement with the theory \cite{Glass} -  details can be found in Appendix \ref{appx:MeanVarFDX}.

\begin{table}[h]
		\centering
	\small
	\begin{tabular}{l|r|r|r|r|r|}
		\cline{2-6}
		& \multicolumn{5}{c|}{\textbf{Time (t)}}                                                                                                              \\ \hline
		\multicolumn{1}{|l|}{$x_0$}    & \multicolumn{1}{c|}{100} & \multicolumn{1}{c|}{1000} & \multicolumn{1}{c|}{5000} & \multicolumn{1}{c|}{10000} & \multicolumn{1}{c|}{20000} \\ \hline
		\multicolumn{1}{|l|}{10}    & 8.35                     & 6.98                      & 2.74                      & 0.00                       & 0.00                       \\ \hline
		\multicolumn{1}{|l|}{100}   & 100.70                   & 110.73                    & 100.92                    & 127.18                     & 149.70                     \\ \hline
		\multicolumn{1}{|l|}{1000}  & 1000.47                  & 996.56                    & 908.82                    & 960.91                     & 1011.71                    \\ \hline
		\multicolumn{1}{|l|}{10000} & 9956.68                  & 9958.66                   & 9752.99                   & 9860.61                    & 9775.20                    \\ \hline
	\end{tabular}
	\caption{Average $x$ value for varying $x_0$ values with $t_0=0$, $\si^2=10$, $\Delta t_i=1$. (10,000 paths)}
	\label{tab:avgxvar10}
\end{table}

\begin{table}[h]
	\centering
	\small
	\begin{tabular}{l|r|r|r|r|r|}
		\cline{2-6}
		& \multicolumn{5}{c|}{\textbf{Time (t)}}                                                                                                              \\ \hline
		\multicolumn{1}{|l|}{$x_0$}    & \multicolumn{1}{c|}{100} & \multicolumn{1}{c|}{1000} & \multicolumn{1}{c|}{5000} & \multicolumn{1}{c|}{10000} & \multicolumn{1}{c|}{20000} \\ \hline
		\multicolumn{1}{|l|}{10}    & 10.07                     & 11.07                      & 10.29                      & 12.68                       & 14.92                       \\ \hline
		\multicolumn{1}{|l|}{100}   & 100.05                   & 99.72                    & 90.96                    & 96.19                     & 99.01                     \\ \hline
		\multicolumn{1}{|l|}{1000}  & 995.67                  & 995.34                    & 975.22                    & 991.77                     & 977.22                    \\ \hline
		\multicolumn{1}{|l|}{10000} & 9990.73                  & 10000.84                   & 10004.78                   & 10049.72                    & 9953.11                    \\ \hline
	\end{tabular}
	\caption{Average $x$ value for varying $x_0$ values with $t_0=0$, $\si^2=1$, $\Delta t_i=1$. (10,000 paths)}
	\label{tab:avgxvar1}
\end{table}

\subsection{Hitting Time Distribution}
\begin{figure}[h]
	\centering
	\subfloat[Feller Sample Path Simulation with $\si^2=1$]{\label{fig:ittingtimea}\includegraphics[width=0.45\linewidth]{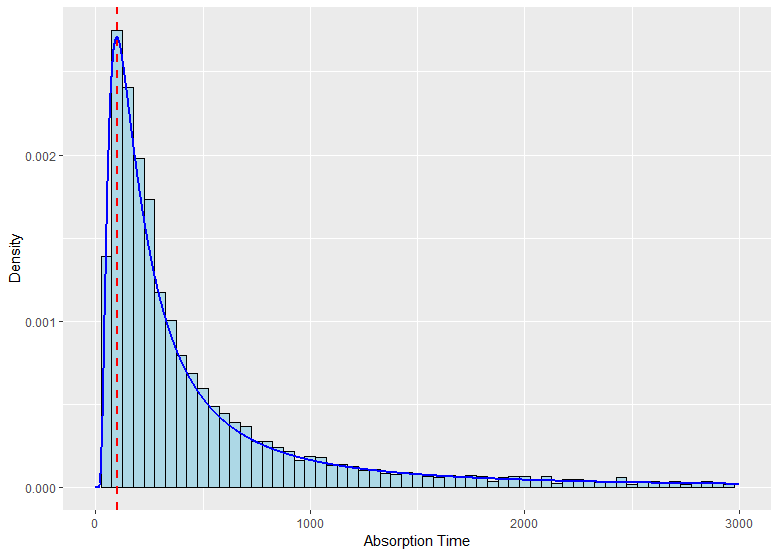}}\qquad
	\subfloat[Feller Sample Path Simulation with $\si^2=10$]{\label{fig:ittingtimeb}\includegraphics[width=0.45\linewidth]{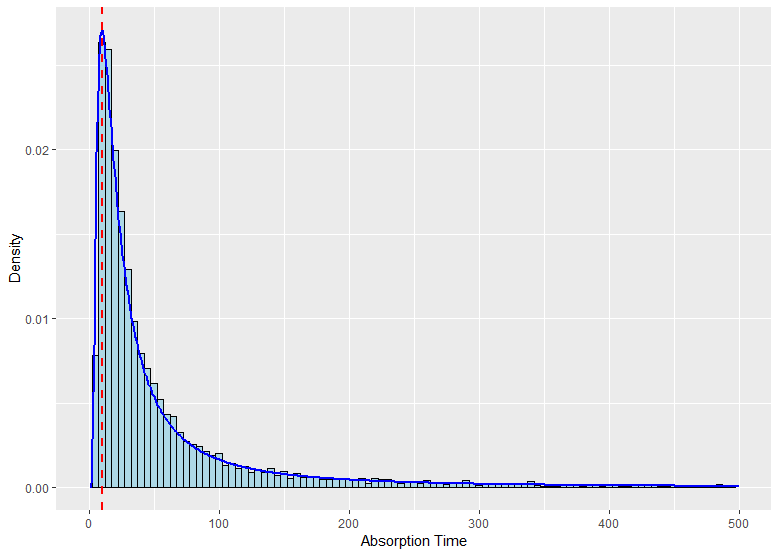}}\\
	\caption{Hitting Time Distribution for Feller walkers with $t_0=0$, $t_n=20000$, $\Delta t_i=1$, $x_0=100$ and different values of $\si^2$.}
	\label{fig:hittingtime}
\end{figure}

From equation (\ref{eqn:absprob}) we defined the absorption probability $F(t)$ and the corresponding survival probability $S(t)$, with the relation $S(t) = 1 - F(t)$. Given the nature of the Feller random walker, we can define these probabilities as:
\bea
S(t) &=& \mathrm{Pr} \left[ x(t) > 0 \right] \label{eqn:S1}  \\
F(t) &=& \mathrm{Pr} \left[ x(t) = 0 \right] \label{eqn:F1}
\eea
We define a new random variable, the hitting time or first-passage time, denoted $T^*$:
\[
T^* = \{ \inf(t) : x(t) = 0 \}
\]
We can then define the hitting (first passage) time density by \cite{Redner}:
\be
\label{eqn:FTPDgen}
f(t) = -\frac{dS(t)}{dt} ,
\ee
where the relation between $f(t)$ and $T^*$ is given by
\[
\mathrm{Probability} \; T^* \in (t, t+dt) = f(t) \; dt
\]
This in turn allows us to express the relations (\ref{eqn:S1}) and (\ref{eqn:F1})  in terns of relationships between $T^*$ and $S(t)$, $F(t)$:
\bea
S(t) &=& \mathrm{Pr} \left[ T^* > t \right] \label{eqn:S2}  \\
F(t) &=& \mathrm{Pr} \left[ T^* < t \right] = \int_{0}^{t} f(\tau) \; d \tau \label{eqn:F2}
\eea
Let us first turn our attention to a Brownian motion originating at $x_o$ described by the SDE $\dot{x}(t) = \si \xi(t)$ with the origin as an absorbing boundary. In this case, the hitting time density is given by \cite{Redner}:
\be
\label{eqn:FTPD2}
f(t) = \frac{x_o}{\sqrt{2} \si t^{3/2}} \exp \left[ \frac{-x_o^2}{2 \si^2 t} \right],
\ee
whose integral over $ t \in [0,\infty)$ is 1, indicating that absorption is certain. Paradoxically however the expected hitting time\footnote{The averaging is with respect to the hitting time density $f(t)$.} $\langle t \rangle = \langle T^* \rangle$ is infinite due to the large time power law behavior $f(t) \sim t^{-3/2}$ \cite{Newman}. Setting $\frac{df}{dt} = 0$ for equation (\ref{eqn:FTPD2}) yields the most likely (typical) hitting time 
\[
t^* = \frac{x_o^2}{3 \si^2}
\]
Now returning to our Feller walker we see an analogous situation to that of the Brownian walker. In this case the hitting time density is:
\be
\label{eqn:FTPD}
f(t) =  \frac{2 x_o}{\si^2 t^2} \exp \left[  \frac{-2 x_o}{\si^2 t}  \right]
\ee
We saw earlier from equation (\ref{eqn:absprob}) that absorption is certain, in terms of the density this is equivalent to the Brownian case with the integral of the hitting density over $t \in [0,\infty)$ equal to 1. As with the Brownian case, due to the large time power law (in this case $f(t) \sim t^{-2}$) the expected hitting time is infinite:
\be
\langle T^* \rangle = \langle t \rangle = \int_{0}^{\infty} t \; f(t) \; dt = \frac{2x_o}{\si^2} \int_{0}^{\infty} \frac{e^{-v}}{v} \; dv
\ee
Similarly setting $\frac{df}{dt} = 0$ for equation (\ref{eqn:FTPD})  yields the typical hitting time for the Feller walker as
\[
t^* = \frac{x_o}{\si^2}
\]
We use simulations to demonstrate the hitting time distributions for a Feller walker with $t_0=0$, $t_n=20000$, $\Delta t_i=1$, $x_0=100$ and $so^2=1,10$ (Figure \ref{fig:hittingtime}). A reference line (shown in red) has been marked on the plots to indicate the theoretical hitting time $t^*$. It can be seen that this appears at the peak of the histograms confirming this result. The theoretical distributions indicated in equation \ref{eqn:FTPD} are overlaid (in blue) as well to illustrate how closely the simulation agrees with the theoretical distribution.

\section{Summary}
In this article, we re-examine the derivation of the probability density of the Feller Diffusion equation. As opposed to the prior approaches \cite{Feller, Feller2, Feller3, Wax}, we use the Fourier transform on the Fokker-Planck equation corresponding to the Feller Diffusion and solve the resulting advection-type equation in Fourier space via the method of characteristics. The solution is easily mapped to the characteristic function of a scaled version of a non-central chi-squared distribution with zero degrees of freedom. As the theory of these distributions is well-known \cite{Siegel}, it is straight forward to interpret the results in terms of survival probabilities and further develop the algorithm for developing sample paths. We carry out such a simulation and i) find the simulated survival probability is very close to the theoretical value and ii) prove the martingale property by estimating the average position of the random walkers using Monte Carlo. We also find the hitting time distribution, including the typical hitting time shows good agreement between the theoretical results and simulation.
For future work we plan to apply the Feller diffusion in the setting of quantitative finance where it has the potential to be used in the following modelling contexts: i) price dynamics of a risky asset or interest rate, ii) derivative pricing on risky underlying assets, iii) bond pricing with risky interest rates,  iv) default probability estimation, v) survival analysis including hazard rate modelling and regression methods.

\section{Acknowledgements}
To the best of our knowledge this work and our approach has not been carried out previously. Our only aim is to add to the body of knowledge in this area. This work has been funded and supported by MIND Analytics \& Management. Contribution to this work as follows: 
\begin{itemize}
\item Motivation for Investigation: RM/LK
\item Fourier Transform + Method of Characteristics Approach: RM
\item Simulation Algorithm Development: RM/LK
\item Simulation, Coding \& Results: PK
\item General Discussion \& Future Direction: RM/LK/PK.
\item Writing: RM/PK/LK
\end{itemize}
The authors declare that there is no conflict of interest regarding the publication of this article.

\appendix

\section{Fourier Transform of Fokker Planck Equation} \label{appx:FTFPE}
In this section we look at deriving equation (\ref{eqn:kFPE}). We begin by re-writing equation (\ref{eqn:FPE}) as
\be
\label{eqn:FPE2}
\frac{\pa}{\pa t} P(x,t) = \frac{1}{2} \si^2 x \frac{\pa^2}{\pa x^2} P(x,t) + \si^2  \frac{\pa}{\pa x} P(x,t) 
\ee
Using the Fourier transform we obtain
\be
\label{eqn:kFPE2}
\frac{\pa}{\pa t} \hat{P}(k,t) = \frac{1}{2} \si^2 \left[ i \frac{\pa}{\pa k} ( - k^2 \hat{P}(k,t) ) \right] + ik \si^2  \frac{\pa}{\pa k} \hat{P}(k,t)
\ee
Note the use of following properties of Fourier transforms in deriving equation (\ref{eqn:kFPE2}):
\bean
\frac{d^nG(x)}{dx^n} \rightarrow (ik)^n \hat{G}(k) \\
x \cdot G(x) \rightarrow i \frac{d}{dk}\hat{G}(k) \\
\eean
Expanding the right hand side of equation  (\ref{eqn:kFPE2}) and simplifying yields the result in equation  (\ref{eqn:kFPE}):
\[
\frac{\pa}{\pa t} \hat{P}(k,t) = \frac{-i \si^2 k^2}{2} \frac{\pa}{\pa k} \hat{P}(k,t)
\]

\section{Solution via Method of Characteristics} \label{appx:MoC}
Let us solve the following PDE for $u(z,t)$ using the method of characteristics:
\bea
\label{eqn:PDE}
\frac{\pa u}{\pa t} + i \alpha z^2 \frac{\pa u}{\pa z} &=& 0 \\
u(z,t_o) &=& \psi(z)
\eea
The method of characteristics postulates a solution of the form $u( z(s), t(s) )$ satisfying
\bea
0      &=& \frac{du}{ds} \label{eqn:MoC} \\
	&=&  \frac{dt}{ds}  \frac{\pa u}{ \pa t} + \frac{dz}{ds}  \frac{\pa u}{ \pa z} \label{eqn:MoCsol}
\eea
Now equation (\ref{eqn:MoC}) implies that along the characteristic curve $(z(s), t(s))$, which is parametrized by $s$, the solution is constant, i.e.
\be 
\label{eqn:char0}
u(s) = u(0) \Rightarrow u( z(s), t(s) ) = u (z(0), t(0)) 
\ee
To be explicit the $(z(s), t(s))$ and $(z(0),t(0))$ lie on the same characteristic curve.
Matching the terms in equation (\ref{eqn:PDE}) with equation (\ref{eqn:MoCsol}) leads to the following system of ordinary differential equations (ODEs) for the characteristics:
\bea
\frac{dt}{ds} &=& 1 \;\;\;t(s=0) = t_o \;\;\;t(s) = t \label{eqn:char1} \\
\frac{dz}{ds} &=& - i \alpha z^2 \;\;\; z(s=0) = z _o \;\;\;z(s) = z \label{eqn:char2} 
\eea
Equation (\ref{eqn:char1}) is solved by integrating as follows:
\be
\label{eqn:char3}
\int_{t_o}^{t} dt' = \int_{0}^{s} ds' \Rightarrow s = t - t_o
\ee
Equation (\ref{eqn:char2}) is solved by integrating as follows:
\be
\label{eqn:char4}
\int_{z_o}^{z} \frac{d \eta}{\eta^2} = -i \alpha \int_{0}^{s} ds' \Rightarrow z_o = \frac{z}{1 + i \alpha z s} = \frac{z}{1 + i \alpha z (t-t_o)} 
\ee
Using equation (\ref{eqn:char0}) along with equations (\ref{eqn:char3}) and (\ref{eqn:char4}) enables us to obtain the solution
\be
\label{eqn:solnU}
u(z,t) = u (z_o, t_o) = \psi(z_o) = \psi \left( \frac{z}{1 + i \alpha z (t-t_o)}  \right) 
\ee
We see that (\ref{eqn:solnU}) satisfies the initial condition $u(z,t_o) = \psi(z)$. For the particular choice of initial condition (as in our problem):
\be
\label{eqn:solnU2}
\psi(z) = e^{-izx_o} \Rightarrow u(z,t) = \exp \left[ \frac{-izx_o}{1 + i \alpha z(t-t_o)} \right]
\ee
We make the identifications to map back to our original problem:
\bean
z	&\rightarrow& k \\
\alpha &\rightarrow& \frac{\si^2}{2} \\
u(z,t) &\rightarrow& \hat{P}(k,t) \\
\psi(z) &\rightarrow& \exp [-i k x_o] \; ,
\eean
It is then straight forward to use (\ref{eqn:solnU2})to arrive at the final answer given by equation (\ref{eqn:soln}):
\[
\hat{P}(k,t) = \exp \left[ \frac{-ikx_o}{1 + \frac{ik \si^2}{2} (t-t_o)} \right]
\]

\section{Mean and Variance of Feller Diffusion Process} \label{appx:MeanVarFDX}
We may use the cumulant generating function $H(k,t)$ to derive the cumulants of the Feller Diffusion, where 
\[
H(k,t) = \ln \varphi (k,t) = \ln \hat{P}(-k,t) = \ln \langle e^{ikx} \rangle  .
\]
In particular the mean is given by
\bea
\langle \; x(t) | \; x(0) = x_o \; \rangle 	&=& -i \frac{\pa }{ \pa k} H(k,t) \bigg\rvert_{k=0} \nn \\
						&=& -i \frac{\pa }{ \pa k} \left[ \frac{-ikx_o}{1 + \frac{ik \si^2}{2} (t-t_o)} \right]  \Bigg\rvert_{k=0} \nn \\
						&=& -i \cdot  ix_o \cdot \left[1 + \frac{ik \si^2}{2} (t-t_o) \right]^{-2}  \bigg\rvert_{k=0} \nn \\
						&=& x_o \label{eq:mart1}
\eea
This matches our intuition as the random walk described by equation (\ref{eqn:SDE}) is drift-less and thus a martingale. The variance of $x(t)$ is given by
\bea
\langle \; \left( x - \langle x \rangle \right)^2  | \; x(0) = x_o \; \rangle 	&=& (-i)^2 \frac{\pa^2 }{ \pa k^2} H(k,t) \bigg\rvert_{k=0} \nn \\
						&=& -i \frac{\pa }{ \pa k} \; x_o \cdot \left[1 + \frac{ik \si^2}{2} (t-t_o) \right]^{-2}   \Bigg\rvert_{k=0} \nn \\
						&=& -i x_o \cdot  i \si^2 (t-t_o) \cdot \left[1 + \frac{ik \si^2}{2} (t-t_o) \right]^{-3}  \Bigg\rvert_{k=0} \nn \\
						&=& \si^2 x_o \cdot (t-t_o) \label{eq:varFDX}
\eea
We use this value to obtain an error-bound on the Monte Carlo estimate of the average position, $\overline{x}(t)$ of a Feller random walker at time $t$ given by equation(\ref{eqn:simmeanpos}).  By Monte Carlo theory \cite{Glass} $\overline{x}(t)$ is distributed (approximately) as a normal with a mean $x_o$ (i.e. true value) and variance $\si^2 x_o t/N$. This enables us to obtain a confidence interval for $x_o$ \cite{Glass}:
\[
\overline{x} (t) \pm \Phi^{-1} ( 1 - \alpha / 2 ) \cdot \sqrt{\frac{\si^2 x_o t}{N}}
\]
Here $\Phi^{-1}$ is the inverse cumulative normal distribution and $\alpha$ is the confidence level\footnote{For the 95 percent  confidence interval $\alpha = 0.05$ and $\Phi^{-1}(0.975) = 1.96$. For the 99 percent confidence interval $\alpha = 0.01$ and $\Phi^{-1}(0.995) = 2.58$}.\\

\noindent To demonstrate this, we include here the results of the Monte Carlo simulations given in Section \ref{sec:MartingaleProperty} with their 95\% confidence intervals in Tables \ref{tab:CIxbarvar10} and \ref{tab:CIxbarvar1}. It is noted that all the simulated results indicate that $\langle x(t)|x(0) = x_o \rangle = x_o$ within the error bounds.

\begin{table}[h]
	\centering
	\small
	\begin{tabular}{|l|l|l|c|c|}
		\hline
		\multicolumn{1}{|c|}{\multirow{2}{*}{$x_0$}} & \multicolumn{1}{c|}{\multirow{2}{*}{Time (t)}} & \multicolumn{1}{c|}{\multirow{2}{*}{$\overline{x}$}} & \multicolumn{2}{c|}{95\% Confidence Interval} \\ \cline{4-5} 
		\multicolumn{1}{|c|}{}                       & \multicolumn{1}{c|}{}                          & \multicolumn{1}{c|}{}                                & Lower Bound           & Upper Bound           \\ \hline
		\multirow{5}{*}{10}                          & 100                                            & 8.35                                                 & 6.39                  & 10.31                 \\ \cline{2-5} 
		& 1000                                           & 6.98                                                 & 0.78                  & 13.18                 \\ \cline{2-5} 
		& 5000                                           & 2.74                                                 & -11.12                & 16.60                 \\ \cline{2-5} 
		& 10000                                          & 0.00                                                 & -19.60                & 19.60                 \\ \cline{2-5} 
		& 20000                                          & 0.00                                                 & -27.72                & 27.72                 \\ \hline
		\multirow{5}{*}{100}                         & 100                                            & 100.70                                               & 94.50                 & 106.90                \\ \cline{2-5} 
		& 1000                                           & 110.73                                               & 91.13                 & 130.33                \\ \cline{2-5} 
		& 5000                                           & 100.92                                               & 57.09                 & 144.75                \\ \cline{2-5} 
		& 10000                                          & 127.18                                               & 65.20                 & 189.16                \\ \cline{2-5} 
		& 20000                                          & 149.70                                               & 62.05                 & 237.35                \\ \hline
		\multirow{5}{*}{1000}                        & 100                                            & 1000.47                                              & 980.87                & 1020.07               \\ \cline{2-5} 
		& 1000                                           & 996.56                                               & 934.58                & 1058.54               \\ \cline{2-5} 
		& 5000                                           & 908.82                                               & 770.23                & 1047.41               \\ \cline{2-5} 
		& 10000                                          & 960.91                                               & 764.91                & 1156.91               \\ \cline{2-5} 
		& 20000                                          & 1011.71                                              & 734.53                & 1288.89               \\ \hline
		\multirow{5}{*}{10000}                       & 100                                            & 9956.68                                              & 9894.70               & 10018.66              \\ \cline{2-5} 
		& 1000                                           & 9958.66                                              & 9762.66               & 10154.66              \\ \cline{2-5} 
		& 5000                                           & 9752.99                                              & 9314.73               & 10191.25              \\ \cline{2-5} 
		& 10000                                          & 9860.61                                              & 9240.81               & 10480.41              \\ \cline{2-5} 
		& 20000                                          & 9775.20                                              & 8898.68               & 10651.72              \\ \hline
	\end{tabular}
	\caption{95\% Confidence Intervals of $\overline{x}$ values given in Table \ref{tab:avgxvar10}}
	\label{tab:CIxbarvar10}
\end{table}

\begin{table}[h]
	\centering
	\small
	\begin{tabular}{|l|l|l|c|c|}
		\hline
		\multicolumn{1}{|c|}{\multirow{2}{*}{$x_0$}} & \multicolumn{1}{c|}{\multirow{2}{*}{Time (t)}} & \multicolumn{1}{c|}{\multirow{2}{*}{$\overline{x}$}} & \multicolumn{2}{c|}{95\% Confidence Interval} \\ \cline{4-5} 
		\multicolumn{1}{|c|}{}                       & \multicolumn{1}{c|}{}                          & \multicolumn{1}{c|}{}                                & Lower Bound           & Upper Bound           \\ \hline
		\multirow{5}{*}{10}                          & 100                                            & 10.07                                                & 8.11                  & 12.03                 \\ \cline{2-5} 
		& 1000                                           & 11.07                                                & 4.87                  & 17.27                 \\ \cline{2-5} 
		& 5000                                           & 10.29                                                & -3.57                 & 24.15                 \\ \cline{2-5} 
		& 10000                                          & 12.68                                                & -6.92                 & 32.28                 \\ \cline{2-5} 
		& 20000                                          & 14.92                                                & -12.80                & 42.64                 \\ \hline
		\multirow{5}{*}{100}                         & 100                                            & 100.05                                               & 93.85                 & 106.25                \\ \cline{2-5} 
		& 1000                                           & 99.72                                                & 80.12                 & 119.32                \\ \cline{2-5} 
		& 5000                                           & 90.96                                                & 47.13                 & 134.79                \\ \cline{2-5} 
		& 10000                                          & 96.19                                                & 34.21                 & 158.17                \\ \cline{2-5} 
		& 20000                                          & 99.01                                                & 11.36                 & 186.66                \\ \hline
		\multirow{5}{*}{1000}                        & 100                                            & 995.67                                               & 976.07                & 1015.27               \\ \cline{2-5} 
		& 1000                                           & 995.34                                               & 933.36                & 1057.32               \\ \cline{2-5} 
		& 5000                                           & 975.22                                               & 836.63                & 1113.81               \\ \cline{2-5} 
		& 10000                                          & 991.77                                               & 795.77                & 1187.77               \\ \cline{2-5} 
		& 20000                                          & 977.22                                               & 700.04                & 1254.40               \\ \hline
		\multirow{5}{*}{10000}                       & 100                                            & 9990.73                                              & 9928.75               & 10052.71              \\ \cline{2-5} 
		& 1000                                           & 10000.84                                             & 9804.84               & 10196.84              \\ \cline{2-5} 
		& 5000                                           & 10004.78                                             & 9566.52               & 10443.04              \\ \cline{2-5} 
		& 10000                                          & 10049.72                                             & 9429.92               & 10669.52              \\ \cline{2-5} 
		& 20000                                          & 9953.11                                              & 9076.59               & 10829.63              \\ \hline
	\end{tabular}
	\caption{95\% Confidence Intervals of $\overline{x}$ values given in Table \ref{tab:avgxvar1}}
	\label{tab:CIxbarvar1}
\end{table}

\section{Zero Degrees of Freedom Non-central Chi-Squared Distribution} \label{appx:0DNCXSq}
A non-central chi-squared distribution with zero degrees of freedom and non-centrality parameter $\lambda$ \cite{Siegel} is given by the density function:
\be \label{eqn:0-Deg}
g(x ; \lambda) = e^{-\lambda / 2} \sum_{j=0}^{\infty} \frac{(\lambda/2)^j}{j!} q(x;2j)
\ee
Here $q(x;\nu)$ is the density function of a central chi-squared distribution with $\nu$ degrees of freedom.
\be
\label{eqn:chisq}
q(x;\nu) = \frac{1}{2^{\nu /2} \Gamma(\nu /2) } x^{\nu / 2 - 1} e^{-x/2}
\ee
Using the convention $q(x;0) = \delta(x)$ we may re-write equation (\ref{eqn:0-Deg}) as
\be \label{eqn:0-Deg1}
g(x ; \lambda) = e^{-\lambda / 2} \left[ \delta(x) + \sum_{j=1}^{\infty} \frac{(\lambda/2)^j}{j!} q(x;2j) \right]
\ee
Equations (\ref{eqn:0-Deg}) and (\ref{eqn:0-Deg1}) describe a compound Poisson mixture of central chi-square distributions with even degrees of freedom \cite{Siegel}. The corresponding characteristic function $\varphi_f (k; \lambda)$ can be derived from (\ref{eqn:0-Deg1}) very easily:
\bea
\varphi_g (k; \lambda)	&=& e^{-\lambda / 2} \left[ 1 + \sum_{j=1}^{\infty} \frac{(\lambda/2)^j}{j!}  \int_{0}^{\infty} q(x;2j) e^{ikx} dx \right] \nn \\
					&=& e^{-\lambda / 2} \left[ 1 + \sum_{j=1}^{\infty} \frac{(\lambda/2)^j}{j!} (1 - 2ik)^{-j} \right] \nn \\
					&=& \exp \left[ \frac{ik \lambda}{1-2ik} \right]
\eea
We have used the fact that the characteristic function of $q(x;2j)$ is $(1 - 2ik)^{-j}$ in the derivation. In addition, we note that we may write the density function in the alternative form \cite{Siegel, Wax}
\be
\label{eqn:MoWax}
g(x;\lambda) = e^{-\lambda / 2} \delta(x) + \frac{1}{2} \sqrt{\frac{\lambda}{x}} e^{-(x+\lambda)/2} I_1 \left(\sqrt{\lambda x} \right)
\ee
Here $I_1(z)$ denotes the modified Bessel function of order 1 which is given by
\be
\label{eqn:I1}
I_1 (z) = \sum_{n=0}^{\infty} \frac{1}{n! (n+1)!}  \left( \frac{z}{2} \right)^{2n+1} = \sum_{n=1}^{\infty} \frac{1}{n! (n-1)!}  \left( \frac{z}{2} \right)^{2n-1}
\ee
Substituting equations (\ref{eqn:chisq}) and (\ref{eqn:I1}) in equation (\ref{eqn:0-Deg1}) yields
\bea
g(x ; \lambda) 	&=& e^{-\lambda / 2} \left[ \delta(x) + \sum_{j=1}^{\infty} \frac{(\lambda/2)^j}{j!} \frac{x^{j-1}}{2^j \Gamma(j)} e^{-x/2} \right] \nn \\
			&=& e^{-\lambda / 2} \delta(x) +  \frac{e^{-(x+\lambda)/2}}{x} \sum_{j=1}^{\infty} \frac{( \lambda x /4)^{j}}{j!(j-1)!} \nn \\
			&=& e^{-\lambda / 2} \delta(x) + \frac{e^{-(x+\lambda)/2}}{x}  \frac{\sqrt{\lambda x}}{2} I_1 \left(\sqrt{\lambda x} \right)
\eea
Re-arranging the last line leads to equation (\ref{eqn:MoWax}). The distribution has a straight forward interpretation where the the probability that random variable takes the value 0 is $e^{-\lambda / 2}$ and the probability that the random variable lies between $a (\neq 0) $ and $b$ is given by the integral of the second term in (\ref{eqn:MoWax}) over the corresponding limits. If $a = 0$ we would include the $\delta$ function term in the integration. Finally, note for the Feller Diffusion process the density function (with $t_o$ set to 0) is given by \cite{Wax}
\be
P(x,t | x_o, 0 ) = e^{-2x_o/ \si^2 t} \cdot \delta(x) + \frac{2}{\si^2 t} \sqrt{\frac{x_o}{x}} \cdot e^ {-2(x+x_o)/ \si^2 t} \cdot I_1 \left(\frac{4 \sqrt{x x_o}}{\si^2 t} \right)
\ee

\section{Rcpp Implementation of Sample Path Simulation} \label{appx:RCPP}
The function CFellerPaths implemented below outputs a list containing time and the set of paths generated. The function allows to specify an initial seed, the number of paths required, a starting and ending time, $\Delta t_i$, the variance and the starting point.

\begin{lstlisting}
#include<Rcpp.h>
#include<math.h>
#include<random>

//[[Rcpp::plugins(cpp11)]]
using namespace Rcpp;

//[[Rcpp::export]]
List CFellerPaths(int seed, int num, int t0, int tn, double dt, double var, 
		double x0){
	int n = ceil(((double(tn - t0))/dt) + 1.0);

	NumericMatrix Paths(n, num);
	static std::mt19937 gen;

	NumericVector t(n);
	t[0] = t0;

	for(int j=0; j<num; j++){
		NumericVector x(n);
		x[0] = x0;

		gen.seed(seed);
		seed += 1;

		for(int i=1; i<n; i++){

			double c = var/(4.0*dt);
			double lambda = float(x[i-1])/c;

			std::poisson_distribution<> pd((lambda/2.0));
			double N = pd(gen);
			if(N==0){
				x[i] = 0.0;
			}else{
				std::chi_squared_distribution<> cd((2*N));
				x[i] = c*cd(gen);
			}
		}

	Paths(_,j) = x;
	}

	for(int i=1; i<n; i++){
		t[i] = t[i-1] + dt;
	}

	return(List::create(Named("t")=t,
	Named("paths")=Paths));

}
\end{lstlisting}

\end{document}